\documentclass[12pt]{article}
\usepackage{amsfonts}
\usepackage{amsmath}
\usepackage{amstext,amsbsy}
\usepackage{epsfig,multicol}
\usepackage{amsthm}
\usepackage{amssymb,latexsym}
\usepackage{color}
\usepackage{eepic}

\allowdisplaybreaks
\newtheorem{theo}{Theorem}[section]

\newtheorem{cor}[theo]{Corollary}

\def\pf{\noindent{\it Proof.} }

\def\qed{\nopagebreak\hfill{\rule{4pt}{7pt}}\medbreak}

\def \P {\mathcal{P}}
\def \Q {\mathcal{Q}}
\def \R {\mathcal{R}}
\def \T {\mathcal{T}}
\def \M {\mathcal{M}}

\def \S {\mathcal{S}}
\def \C {\mathcal{C}}

\begin{document}

\begin{center}
{\large\bf Partitions and Partial Matchings
Avoiding Neighbor Patterns }\\
\end{center}
\begin{center}
William Y. C. Chen$^{1}$, Neil J. Y. Fan$^{2}$ and
Alina F. Y. Zhao$^{3}$\\[6pt]
Center for Combinatorics, LPMC-TJKLC\\
Nankai University, Tianjin 300071, P. R. China\\[5pt]
$^{1}${chen@nankai.edu.cn}, $^{2}${fjy@cfc.nankai.edu.cn},
$^{3}${zfeiyan@cfc.nankai.edu.cn}
\end{center}

\vskip 3mm \noindent \textbf{Abstract.}
  We obtain the generating functions for partial matchings avoiding neighbor alignments and for partial matchings avoiding neighbor alignments and left nestings.
 We show that there is a bijection between partial
  matchings avoiding three neighbor patterns
  (neighbor alignments, left nestings and right nestings)
  and
  set partitions avoiding right nestings via an intermediate
  structure of integer compositions.
  Such integer compositions
  are known to be in one-to-one correspondence
  with self-modified ascent sequences
  or $3\bar{1}52\bar{4}$-avoiding permutations, as shown
  by Bousquet-M\'elou, Claesson, Dukes and Kitaev.



%%%%%%%%%%%%%%%%%%%%%%%%%%%%%%
%  1. Introduction           %
%%%%%%%%%%%%%%%%%%%%%%%%%%%%%%

\section{Introduction } \label{S:introduction}

This paper is concerned with the enumeration of
partial matchings and set partitions that avoid certain  neighbor
patterns. Recall that a partition $\pi$ of
$[n]=\{1,2,\ldots,n\}$
can be represented as a diagram with
 vertices drawn on a horizontal line
  in increasing order. For a block $B$ of $\pi$,
  we write the elements of $B$ in increasing order.
  Suppose that $B=\{i_1, i_2, \ldots, i_k\}$.
  Then we draw an arc from $i_1$ to $i_2$,
  an arc from $i_2$ to $i_3$, and so on.
   Such a diagram is called the linear representation
    of $\pi$. If $(i,j)$ is an arc in the diagram
    of $\pi$, we call $i$ a left-hand endpoint,
    and call $j$ a right-hand endpoint.

A partial matching is a partition for
which each block contains at
most two elements. A partial matching is also
 called a poor
partition by Klazar  \cite{Kla98},
see also   \cite{Chen05}, and
it can be viewed as an involution on a set.
 A partition for which each block
contains exactly two elements is called a perfect matching.

 Perfect matchings avoiding
 certain patterns have been studied in
\cite{Chen06,Chen07,CL,M,JLMY,Kla1,Kla2,Kla3,Stein}.
Bousquet-M\'elou, Claesson, Dukes and Kitaev \cite{BCDK}
investigated perfect matchings avoiding left nestings
 and right
nestings, and found bijections with other
combinatorial objects such
as $(2+2)$-free posets.
   Claesson and Linusson \cite{CL}
established a correspondence between permutations
and perfect matchings avoiding left nestings.

A nesting
of a partition $\pi$ is formed by two
arcs $(i_1,j_1)$ and $(i_2,
j_2)$
  in the linear representation
   such that $i_1<i_2<j_2<j_1$.
    If we further require that $i_1+1=i_2$, then this nesting is called a
left nesting. Similarly, one can define right nestings, as well as
left crossings and right crossings.  We say that
 $k$ arcs $(i_1,
j_1),(i_2, j_2),\ldots,(i_k, j_k)$ form a $k$-crossing if $i_1 <
i_2<\cdots<i_k< j_1 < j_2<\cdots<j_k$. An alignment of a partition
$\pi$ is formed by two arcs $(i_1,j_1)$ and $(i_2, j_2)$ such that
$i_1<j_1<i_2<j_2$.

\begin{figure}[h,t]\label{fig1}
\begin{center}
\begin{picture}(410,20)
\setlength{\unitlength}{0.5mm}

\put(0,0){\circle*{1.5}}\put(10,0){\circle*{1.5}}\put(15,-2){$\cdots$}
\put(28,0){\circle*{1.5}}\put(32,-2){$\cdots$}\put(45,0){\circle*{1.5}}
\put(-2,-9){\footnotesize$i$}\put(5,-9){\footnotesize$i+1$}\put(25,-9){\footnotesize$j_1$}
\put(42,-9){\footnotesize$j_2$}
\qbezier[1000](0,0)(14,20)(28,0)\qbezier[1000](10,0)(27.5,25)(45,0)

\put(60,0){\circle*{1.5}}\put(70,0){\circle*{1.5}}\put(75,-2){$\cdots$}
\put(88,0){\circle*{1.5}}\put(92,-2){$\cdots$}\put(105,0){\circle*{1.5}}
\put(58,-9){\small$i$}\put(65,-9){\footnotesize$i+1$}\put(85,-9){\footnotesize$j_1$}
\put(102,-9){\footnotesize$j_2$}
\qbezier[1000](60,0)(82.5,25)(105,0)\qbezier[1000](70,0)(79,10)(88,0)

\put(120,0){\circle*{1.5}}\put(125,-2){$\cdots$}\put(138,0){\circle*{1.5}}
\put(142,-2){$\cdots$}\put(155,0){\circle*{1.5}}\put(165,0){\circle*{1.5}}
\put(118,-9){\footnotesize$i_1$}\put(136,-9){\footnotesize$i_2$}\put(153,-9){\footnotesize$j$}
\put(159,-9){\footnotesize$j+1$}
\qbezier[1000](120,0)(137.5,25)(155,0)\qbezier[1000](138,0)(151.5,20)(165,0)

\put(180,0){\circle*{1.5}}\put(185,-2){$\cdots$}\put(198,0){\circle*{1.5}}
\put(202,-2){$\cdots$}\put(215,0){\circle*{1.5}}\put(225,0){\circle*{1.5}}
\put(178,-9){\footnotesize$i_1$}\put(196,-9){\footnotesize$i_2$}\put(213,-9){\footnotesize$j$}
\put(219,-9){\footnotesize$j+1$}
\qbezier[1000](180,0)(202.5,25)(225,0)\qbezier[1000](198,0)(206.5,10)(215,0)

\put(240,0){\circle*{1.5}}
\put(245,-2){$\cdots$}\put(258,0){\circle*{1.5}}\put(268,0){\circle*{1.5}}\put(273,-2){$\cdots$}
\put(286,0){\circle*{1.5}}
\put(238,-9){\footnotesize$i_1$}\put(255,-9){\footnotesize$j_1$}\put(263,-9){\footnotesize$j_1+1$}
\put(284,-9){\footnotesize$j_2$}
\qbezier[1000](240,0)(249,20)(258,0)\qbezier[1000](268,0)(275,20)(286,0)

\end{picture}
\end{center}
\caption{Left crossing, left nesting, right crossing, right nesting
and neighbor alignment.}
\end{figure}
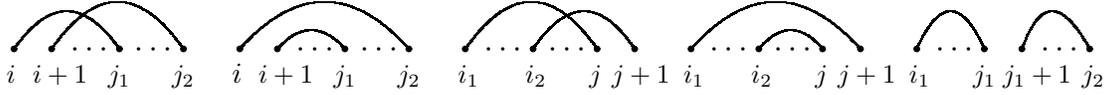

In this paper, we define a neighbor alignment  as
an alignment consisting of two arcs $(i_1,j_1)$
 and $(i_2,j_2)$ such
that $j_1+1=i_2$. The aforementioned patterns with neighbor
constraints are called neighbor patterns. Left nestings and right
nestings were introduced by Stoimenow \cite{Sto} in the study of
regular linearized chord diagrams, and were further explored in
\cite{BCDK,CL,CDK}. An illustration of   neighbor patterns is given
in Figure \ref{fig1}.

Our main results are the generating functions for three classes of partial matchings avoiding neighbor patterns, which are denoted by
  $\P(n),\Q(n),\R(n)$, respectively.
Denote the set of partial matchings of $[n]$ by $\M(n)$. The  set of partial matchings in $\M(n)$ with no neighbor alignments is denoted by $\P(n)$, and
the set of  partial matchings in $\P(n)$ with $k$ arcs is
denoted by $\P(n,k)$. The set of
 partial matchings in $\P(n)$ with no left
nestings is denoted by $\Q(n)$, and the
 set of partial matchings in $\Q(n)$ with $k$ arcs is denoted by  $\Q(n,k)$.
Moreover, the set of partial matchings in $\Q(n)$ with no right nestings is denoted by $\R(n)$,
  and the set of partial matchings in  $\R(n)$ with $k$ arcs is denoted by $\R(n,k)$. For $0\leq k \leq \lfloor  n/2\rfloor$, we set
    \[ P(n,k)=|\P(n,k)|, \quad Q(n,k)=|\Q(n,k)|, \quad R(n,k)=|\R(n,k)|.\]
Denote the set of partitions of $[n]$ by $\S(n)$ and denote the set of
partitions in $\S(n)$ with $k$ blocks  by $\S(n,k)$. The
set of partitions in $\S(n)$ with no right nestings is denoted by
$\T(n)$, and the set of partitions in $\T(n)$ with $k$ arcs is
denoted by $\T(n,k)$. For $0\leq k \leq n-1$, we set
$T(n,k)=|\T(n,k)|$.

We obtain the following generating function
formulas for the numbers $P(n,k)$ and $Q(n,k)$.

\begin{theo}\label{t1}

\begin{equation}\label{gf1}
\sum_{n\geq 1}\sum_{k=0}^{\lfloor \frac{n}{2}\rfloor}
P(n,k)x^ny^k=\sum_{n\geq 1}\prod_{k=1}^n (1+kxy)x^n.
\end{equation}
\end{theo}

\begin{theo}\label{t2}
\begin{align}\label{gf2}
\sum_{n\geq 1}\sum_{k=0}^{\lfloor \frac{n-1}{2}\rfloor}
Q(n-1,k)x^ny^k=\sum_{n\geq 1} \frac{x^n}{\prod_{k=1}^n(1-kx^2y)}.
\end{align}
\end{theo}

It is clear that when $y=1$,  the right-hand
side of \eqref{gf1} reduces to
\[\sum_{n\geq 1}\prod_{k=1}^n (1+kx)x^n\]
which is the generating function of the sequence $A124380$ in OEIS
\cite{slo}, whose first few entries are
\[1, 2, 4, 9, 22, 57, 157, 453, 1368, 4290, \ldots.\]
It seems that no combinatorial interpretations of this sequence are
known.  Thus Theorem \ref{t1} can be considered as a combinatorial
interpretation of  the above generating function.

 Meanwhile, when $y=1$
the right-hand side of \eqref{gf2} reduces to
\[\sum_{n\geq 1} \frac{x^n}{\prod_{k=1}^n(1-kx^2)}\]
which is the generating function of the sequence $A024428$ in OEIS \cite{slo}, whose first
few entries are
\[1,1, 2, 4, 8, 18, 42, 102, 260, 684, 1860, \ldots.\]
This sequence can be expressed  in terms of Stirling numbers of the
second kind. So Theorem \ref{t2} can be considered as another
combinatorial interpretation  of the above generating function.

We derive the generating function for the numbers $R(n,k)$ by
establishing a connection with compositions of the integer $n-k$
into ${k+1 \choose 2}$ components. Moreover, we show that there is a
correspondence between the set $\R(n,k)$ and the set $\T(n-k+1,k)$.
Hence by Theorem \ref{t3} we obtain the generating function for
$T(n,k)$ as stated in Theorem \ref{t4}.
Furthermore, it turns out that
this generating
 function coincides with the generating function
 for the number of
 self-modified ascent sequences of length
 $n$ with largest element $k-1$
 or $3\bar{1}52\bar{4}$-avoiding
 permutations having $k$ right-to-left minima, as
 derived by Bousquet-M\'elou, Claesson,
 Dukes and Kitaev \cite{BCDK}.
\begin{theo}\label{t3}
\begin{align}\label{gf3}\sum_{n\geq
1}\sum_{k=0}^{n-1}R(n+k-1,k)x^ny^{k}=\sum_{n\geq 1}
\frac{x^n}{(1-xy)^{{n+1 \choose 2}}}.
\end{align}

\end{theo}

\begin{theo}\label{t4}
\begin{align}\label{gf4}
\sum_{n\geq1}\sum_{k=0}^{n-1}T(n,k)x^ny^{k}=\sum_{n\geq 1}
\frac{x^n}{(1-xy)^{{n+1 \choose 2}}}.
\end{align}
\end{theo}

This paper is structured as follows. In Section 2, we give a proof
of Theorem 1.1 by deriving a recurrence relation of $P(n,k)$.
Section 3 gives a proof of Theorem 1.2 by establishing  a
correspondence between $\S(n-k,n-2k)$ and $\Q(n-1,k)$. In Section 4, we give
a bijection between $\R(2n-k-1,n-k)$ and the set of compositions of
 $n-k$ into ${k+1 \choose 2}$ components, which leads to the generating function in Theorem
1.3. In Section 5 we present a proof of Theorem 1.4
by constructing a
correspondence between the set $\R(n,k)$
and the set $\T(n-k+1,k)$.

\section{Neighbor alignments}

In this section, we give a proof of the generating function formula for the number of
partial matchings avoiding neighbor alignments. Recall that a
singleton of a partial matching or a set partition is the only
element in a block, which corresponds to an isolated vertex in its
diagram representation. For a block with at least two elements, the minimum element
is called an origin, and  the maximum element  is called a
destination, and an element in between, if any, is called a
transient. An origin and a destination are also called an opener and
a closer respectively by some authors.
 We first give a recurrence relation of  $P(n,k)$.

\begin{theo}\label{p1}
For $n\geq 3$, and $1\leq k \leq n/2$, we have
\begin{equation} \label{rec1}
P(n,k)=P(n-1,k)+(n-k)P(n-2,k-1),
\end{equation}
with initial values $P(1,0)=1, P(2,0)=1, P(2,1)=1.$
\end{theo}

\pf  It is clear that the number of partial
  matchings in $\P(n,k)$ such that the element $1$
  is a singleton equals $P(n-1, k)$. So it suffices
  to show that the number of partial matchings in which $1$
   is not a singleton equals $(n-k)P(n-2, k-1)$.
   For a partial matching $M \in \P(n,k)$ in which $1$ is not a
singleton,
   after deleting the arc with left-hand endpoint $1$, we are led to a
   partial matching in $\P(n-2,k-1)$.

   Conversely, given a partial matching
   $M \in \P(n-2,k-1)$ with
   $n-2$ vertices, in order to
get a   partial matching with $k$ arcs,
we can add an arc into $M$ by placing the
    left-hand endpoint
  before the first vertex of $M$ and inserting
right-hand endpoint at some position of $M$. Clearly,
there are $n-1$ possible positions to insert the
right-hand endpoint of the new arc. By
abuse of language, if no confusion arises we do not distinguish a
partial matching $M$ from its diagram representation. To ensure
that the insertion will not cause any neighbor alignments,
 we should not allow
 the right-hand endpoint of the
inserted arc  to be placed before any origin
of $M$. Since there are $k-1$ arcs in  $M$,
thus there are $k-1$
positions that are forbidden.
That is to say, we  have exactly
$(n-1)-(k-1)=n-k$ choices for the position of the right-hand
endpoint of the inserted arc.
 After relabeling, we get a partial
matching in $P(n,k)$. This completes the proof.
 \qed

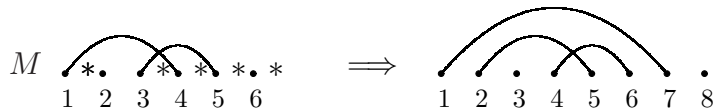
\begin{figure}[h,t]\label{fig2}
\begin{center}
\begin{picture}(235,20)
\setlength{\unitlength}{0.5mm}

\put(-15,0){$M$}

\multiput(0,0)(10,0){6}{\circle*{1.5}}
\put(-1,-9){\footnotesize1}\put(9,-9){\footnotesize2}\put(19,-9){\footnotesize3}\put(29,-9){\footnotesize4}
\put(39,-9){\footnotesize5}\put(49,-9){\footnotesize6}

\qbezier[1000](0,0)(15,20)(30,0) \qbezier[1000](20,0)(30,15)(40,0)

\put(4,-1){$*$}\put(4,-1){$*$}\put(24,-1){$*$}
\put(34,-1){$*$}\put(44,-1){$*$}\put(54,-1){$*$}

\put(75,0){$\Longrightarrow $}

\multiput(100,0)(10,0){8}{\circle*{1.5}}

\put(99,-9){\footnotesize1}\put(109,-9){\footnotesize2}\put(119,-9){\footnotesize3}
\put(129,-9){\footnotesize4}\put(139,-9){\footnotesize5}
\put(149,-9){\footnotesize6}\put(159,-9){\footnotesize7}\put(169,-9){\footnotesize8}
\qbezier[1000](100,0)(130,35)(160,0)
\qbezier[1000](110,0)(125,20)(140,0)
\qbezier[1000](130,0)(140,15)(150,0)
\end{picture}
\vspace*{5mm} \caption{Possible positions for inserting an arc.}
\end{center}
\end{figure}

As an example, let us consider a partial matching
$M=\{\{1,4\},\{2\},\{3,5\},\{6\}\}\in \P(6,2)$. The possible
positions for inserting an arc are marked by the symbol $*$ in
Figure \ref{fig2}.¡¡Note that the positions before the vertices $1$
and $3$ are forbidden. The right-hand endpoint of the inserted arc
is between  the vertices $5$ and $6$.

%%%%%%%%%%%%%%%%%%%%%%%%%%%%%%%%%%%%%%%%%%%%%%%%%%%%%%%%%
{\noindent \it Proof of Theorem 1.1.} Let $f(n,k)$ denote the coefficient of $x^{k}y^k$ in the expansion of
 \[ \prod_{i=1}^{n-k} (1+ixy). \] It is easily verified that
\[f(n,k)=f(n-1,k)+(n-k)f(n-2,k-1),\]
 with initial values
\[f(1,0)=1,\quad
f(2,0)=1,\quad f(2,1)=1.\] Consequently,  $P(n,k)$ and $f(n,k)$ have
the same recurrence relation and the same initial values, so they
are equal. This completes the proof. \qed

To conclude this section, we give a recurrence
relation of the generating function of $P(n,k)$.
Let \[ f_n(y)=\sum_{k=0}^{\lfloor\frac{n}{2}\rfloor} P(n,k)y^k.\]

\begin{cor}
For $n\geq 3$, we
have
\begin{equation}
 f_n(y)=f_{n-1}(y)+(n-1)yf_{n-2}(y)-y^2f_{n-2}'(y).
\end{equation}
\end{cor}

\section{Neighbor alignments and left nestings}

This section is concerned with the
generating function for partial matchings
avoiding neighbor alignments and left nestings.
 More precisely, we
establish a  bijection between set partitions and partial matchings
avoiding neighbor alignments and left nestings. As a consequence, we
obtain the generating function in Theorem \ref{t2}.

\begin{theo}\label{3-1}
There exists a bijection between the set $\S(n-k,n-2k)$ and the set
$\Q(n-1,k)$. Moreover, this bijection transforms the number of
transients of a partition to the number of left crossings of a
partial matching.
\end{theo}

\pf Let $\pi\in \S(n-k,n-2k)$ be a partition
 of $[n-k]$ with $k$
arcs, we wish to add $k-1$ vertices  to $\pi$
in order to form a partial
matching $\alpha(\pi)$ in $\Q(n-1,k)$,
 that is, a partial matching
on $[n-1]$ avoiding neighbor alignments and left nestings.
First, we
add a vertex before each origin,
 except for the first origin, and
relabel the vertices in the natural order.
 Let the resulting partition be
denoted by $\sigma$.

\begin{figure}[h,t]
\begin{center}
\begin{picture}(200,20)
\setlength{\unitlength}{1mm} \multiput(0,0)(10,0){3}{\circle*{.8}}
\put(2.5,-1){$\cdots$} \qbezier[1000](0,0)(5,10)(10,0)
\qbezier[1000](10,0)(15,10)(20,0)
\put(12.5,-1){$\cdots$}\put(-1,-4){\small$i$}\put(9,-4){\footnotesize$j$}\put(19,-4){\footnotesize$k$}

 \put(26,1){$\Longrightarrow$}
\multiput(40,0)(5,0){2}{\circle*{.8}}\multiput(55,0)(10,0){2}{\circle*{.8}}
\qbezier[1000](40,0)(47.5,10)(55,0)
\qbezier[1000](45,0)(55,10)(65,0) \put(47.5,-1){$\cdots$}
\put(57.5,-1){$\cdots$}
\put(39,-4){\footnotesize$i$}\put(42.5,-4){\footnotesize$i+1$}
\put(52.5,-4){\footnotesize$j+1$}\put(62.5,-4){\footnotesize$k+1$}

\end{picture}
\vspace*{.3cm}
 \caption{Change a 2-path into a left crossing.}\label{lc}
\end{center}
\end{figure}
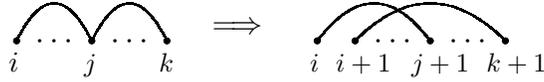

To construct a partial matching in $\Q(n-1,k)$ from the
partition $\sigma$, we shall use the operation
 of changing a 2-path
to a left crossing, see Figure \ref{lc} for an illustration.
To be
more specific, a 2-path means two arcs $(i,j)$ and $(j,k)$ with
$i<j<k$.

It should be emphasized that at each step we have a unique choice of
a $2$-path to implement the switching operation. More precisely, we
always try to find a $2$-path consisting of $(i,j)$ and $(j,k)$
 such that $j$ is
minimal.  We add a vertex $i+1$ immediately after the vertex $i$,
transform this 2-path into a left crossing $(i,j+1),(i+1,k+1)$ and
 relabel those vertices that are greater than $i$ (except for
$j,k$) so that the relabeled vertex set becomes a set of the
 first
consecutive natural numbers.
 In other words, a vertex is relabeled
if it is increased by one. Using this operation to a path
corresponding to a block $B$ with $r+1$ elements
in a partition
$\pi$, we get a left $r$-crossing, which is
an $r$-crossing with
consecutive left-hand endpoints.
 Let $\alpha(\pi)$ denote the
resulting partial matching.

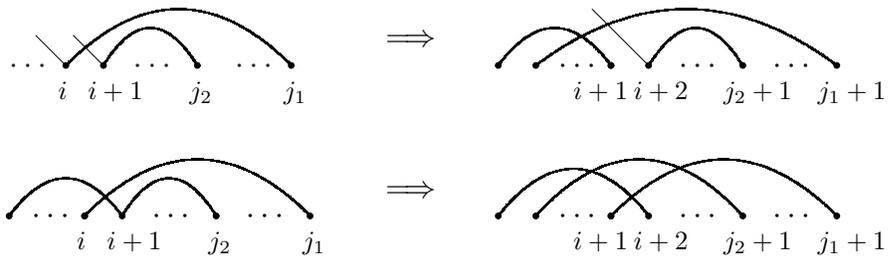
\begin{figure}[h,t]
\begin{center}
\begin{picture}(300,80)
\setlength{\unitlength}{0.5mm}

\put(0,38){$\cdots$}\put(33,38){$\cdots$}\put(60,38){$\cdots$}
\put(15,40){\circle*{1.5}}\put(25,40){\circle*{1.5}}
\put(50,40){\circle*{1.5}} \put(75,40){\circle*{1.5}}
\put(13,31){\footnotesize$i$}\put(21,31){\footnotesize$i+1$}
\put(48,31){\footnotesize$j_2$}\put(73,31){\footnotesize$j_1$}
\qbezier[1000](15,40)(45,70)(75,40)\qbezier[1000](25,40)(37.5,60)(50,40)
\put(15,40){$\line(-1,1){8}$}\put(25,40){$\line(-1,1){8}$}

\put(100,45){$\Longrightarrow $}

\put(130,40){\circle*{1.5}}\put(140,40){\circle*{1.5}}
\put(146,38){$\cdots$}\put(160,40){\circle*{1.5}}
\put(178,38){$\cdots$}\put(203,38){$\cdots$}
\put(160,40){\circle*{1.5}}\put(170,40){\circle*{1.5}}
\put(195,40){\circle*{1.5}} \put(220,40){\circle*{1.5}}
\put(150,31){\footnotesize$i+1$}\put(166,31){\footnotesize$i+2$}
\put(190,31){\footnotesize$j_2+1$}\put(215,31){\footnotesize$j_1+1$}

\qbezier[1000](130,40)(145,60)(160,40)\qbezier[1000](140,40)(180,70)(220,40)
\qbezier[1000](170,40)(182.5,60)(195,40)
\put(170,40){$\line(-1,1){15}$}

%%%%%%%%%%%%%

\put(0,0){\circle*{1.5}}\put(6,-2){$\cdots$}\put(20,0){\circle*{1.5}}
\put(38,-2){$\cdots$}\put(63,-2){$\cdots$}
\put(20,0){\circle*{1.5}}\put(30,0){\circle*{1.5}}
\put(55,0){\circle*{1.5}} \put(80,0){\circle*{1.5}}
\put(18,-9){\footnotesize$i$}\put(26,-9){{\footnotesize{$i+1$}}}
\put(53,-9){\footnotesize$j_2$}\put(78,-9){\footnotesize$j_1$}
\qbezier[1000](0,0)(15,20)(30,0)\qbezier[1000](20,0)(50,30)(80,0)
\qbezier[1000](30,0)(42.5,20)(55,0)

\put(100,5){$\Longrightarrow $}

\put(130,0){\circle*{1.5}}\put(140,0){\circle*{1.5}}
\put(146,-2){$\cdots$}\put(160,0){\circle*{1.5}}
\put(178,-2){$\cdots$}\put(203,-2){$\cdots$}
\put(160,0){\circle*{1.5}}\put(170,0){\circle*{1.5}}
\put(195,0){\circle*{1.5}} \put(220,0){\circle*{1.5}}
\put(150,-9){\footnotesize$i+1$}\put(166,-9){\footnotesize$i+2$}
\put(190,-9){\footnotesize$j_2+1$}\put(215,-9){\footnotesize$j_1+1$}

\qbezier[1000](130,0)(150,25)(170,0)
\qbezier[1000](140,0)(167.5,30)(195,0)
\qbezier[1000](160,0)(190,30)(220,0)

\end{picture}
\vspace*{.5cm}
 \caption{The two cases for the vertex $i$.}\label{fig3}
\end{center}
\end{figure}

We claim that there are no left nestings and neighbor alignments in
$\alpha(\pi)$. Recall that after the first step,
a possible left
nesting in $\sigma$
consisting of arcs $(i,j_1)$ and $(i+1,j_2)$
with $j_1>j_2$ can  occur only when $i+1$ is a transient.
Clearly,
 $i$ is either a transient or an origin. After
changing the 2-path for which $i$ is a transient into a left
  crossing and the 2-path for which $i+1$ is a transient into a left
  crossing, we see that
  the left
nesting $(i,j_1),(i+1,j_2)$ disappears.
 This operation is
illustrated in Figure \ref{fig3}.
A possible neighbor alignment
consisting of arcs $(i,j)$ and $(j+1,k)$
in $\sigma$ can occur only
when $j+1$ is a transient.
 After changing the 2-path for
which $j+1$ is a transient into
a left crossing, the arc $(j+1,k)$
becomes $(j',k+1)$ with $j'<j$, thus the neighbor alignment
disappears. Hence the claim is proved.

It remains  to show that there are $n-1$ vertices
 in $\alpha(\pi)$.
Adding a vertex immediately before an origin or transforming
 a
2-path into a left crossing will result
 in an increase of the number
of vertices by one. Since the left-hand endpoint of an arc is either
an origin or a transient,  for a partition with $k$ arcs,
 after these two operations there are
  a total number of $k-1$ vertices
added. This implies that  $\alpha(\pi) \in \Q(n-1,k)$.

Conversely, given a partial matching in $\Q(n-1,k)$, in order to
obtain a partition in $\S(n-k,n-2k)$,
we must delete $k-1$ vertices.
First, we change   left crossings to 2-paths  from
right to left until there are no more left crossings.
More
precisely, suppose that the two arcs $(i,j_1)$
and $(i+1,j_2)$ form
a left crossing which is the rightmost one in the sense that
 $i$ is the
largest among all the first origins of left crossings.
For such a
left crossing, delete the vertex $i+1$, and
 change the two arcs
$(i,j_1)$ and $(i+1,j_2)$ to $(i,j_1-1)$
 and $(j_1-1,j_2-1)$. Then,
relabel the vertices, that is,
reduce the labels of vertices  larger
than $i$ (except $j,k$) by 1. After we have eliminated
 all the left
crossings, we further delete the singleton
 immediately before each
origin, if there is any,  except for the singleton
 immediately before the
first origin. Finally, we relabel the vertices  by the natural
order.

Since there are neither neighbor alignments
nor left nestings in a
partial matching in $\Q(n-1,k)$ with $k$ arcs,
 for any arc $(i,j)$, we have only two possibilities:
 (1) The vertex $i-1$ is a singleton.
 (2) There exists a vertex $k$ such that
     $(i-1, k)$ and $(i,j)$ form a left crossing.
 Transforming  a left crossing into a 2-path or deleting a vertex immediately before an
origin will result
 in a decrease of the number
of vertices by one. Thus after these two operations there are a
total number of $k-1$ vertices deleted, and so we are led to a
partition of $[n-k]$ with $k$ arcs. It is easily seen that the
number of transients  of $M$ equals the number of left crossings of
$\alpha(M)$. This completes the proof.
 \qed

\begin{figure}[h,t]
\begin{center}
\begin{picture}(300,160)
\setlength{\unitlength}{0.5mm}

\put(20,0){$\Longrightarrow $}
\multiput(50,0)(10,0){12}{\circle*{1.5}}
\put(49,-9){\footnotesize1}\put(59,-9){\footnotesize2}\put(69,-9){\footnotesize3}\put(79,-9){\footnotesize4}\put(89,-9){\footnotesize5}
\put(99,-9){\footnotesize6}\put(109,-9){\footnotesize7}
\put(119,-9){\footnotesize8}\put(129,-9){\footnotesize9}\put(137,-9){\footnotesize10}\put(147,-9){\footnotesize11}\put(157,-9){\footnotesize12}
\qbezier[1000](50,0)(85,40)(120,0)
\qbezier[1000](70,0)(85,25)(100,0)\qbezier[1000](80,0)(95,25)(110,0)
\qbezier[1000](90,0)(120,35)(150,0)\qbezier[1000](140,0)(150,20)(160,0)

\put(20,50){$\Longrightarrow $}
\multiput(50,50)(10,0){11}{\circle*{1.5}}
\put(49,41){\footnotesize1}\put(59,41){\footnotesize2}\put(69,41){\footnotesize3}\put(79,41){\footnotesize4}\put(89,41){\footnotesize5}\put(99,41){\footnotesize6}
\put(109,41){\footnotesize7}\put(119,41){\footnotesize8}\put(129,41){\footnotesize9}\put(138,41){\footnotesize10}\put(148,41){\footnotesize11}
\qbezier[1000](50,50)(80,85)(110,50)\qbezier[1000](70,50)(80,70)(90,50)
\qbezier[1000](80,50)(90,70)(100,50)\qbezier[1000](100,50)(120,80)(140,50)
\qbezier[1000](130,50)(140,70)(150,50)

\multiput(0,100)(10,0){8}{\circle*{1.5}}
\put(-1,91){\footnotesize1}\put(9,91){\footnotesize2}\put(19,91){\footnotesize3}\put(29,91){\footnotesize4}\put(39,91){\footnotesize5}\put(49,91){\footnotesize6}
\put(59,91){\footnotesize7}\put(69,91){\footnotesize8}
\qbezier[1000](0,100)(20,130)(40,100)\qbezier[1000](10,100)(15,120)(20,100)
\qbezier[1000](30,100)(25,120)(20,100)\qbezier[1000](30,100)(45,125)(60,100)
\qbezier[1000](50,100)(60,120)(70,100)

\put(93,100){$\Longrightarrow $}

\multiput(130,100)(10,0){10}{\circle*{1.5}}
\put(129,91){\footnotesize1}\put(139,91){\footnotesize2}\put(149,91){\footnotesize3}\put(159,91){\footnotesize4}\put(169,91){\footnotesize5}
\put(179,91){\footnotesize6}\put(189,91){\footnotesize7}\put(199,91){\footnotesize8}\put(209,91){\footnotesize9}\put(218,91){\footnotesize10}
\qbezier[1000](130,100)(155,130)(180,100)
\qbezier[1000](150,100)(155,120)(160,100)\qbezier[1000](160,100)(165,120)(170,100)
\qbezier[1000](170,100)(190,130)(210,100)\qbezier[1000](200,100)(210,120)(220,100)

\end{picture}
\vspace*{5mm} \caption{The bijection $\alpha$.}\label{fig4}
\end{center}
\end{figure}
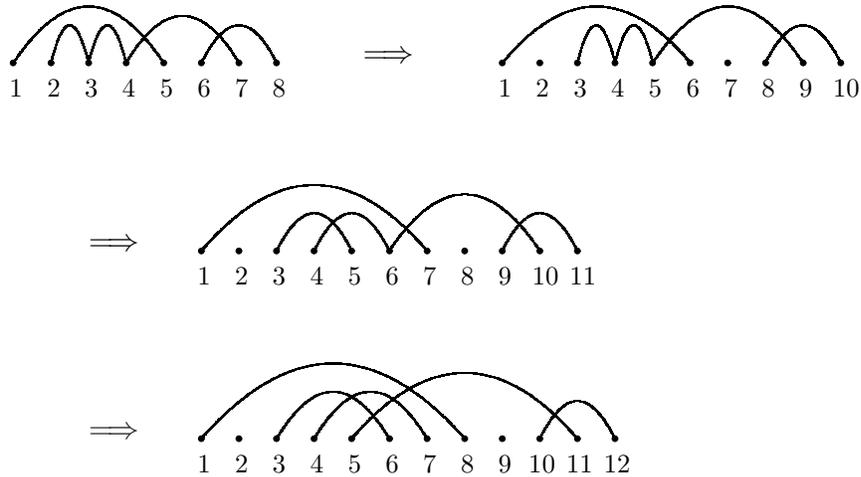

For example, let $\pi=\{\{1,5\},\{2,3,4,7\},\{6,8\}\}\in \S(8,3)$.
There are $5$ arcs in $\pi$, so we must
 add $4$ vertices   in order to get
a partial matching in $\Q(12,5)$.
We first add a singleton before
the arc $(2,3)$ and  a singleton before the arc $(6,8)$.
Then change
the two 2-paths into left crossings
 from left to right.
An illustration of the above procedure is shown in Figure
\ref{fig4}.

Let $g_n(y)$ be the generating function for the numbers $Q(n-1,k)$.
From Theorem \ref{3-1} we see that
\[g_n(y) =\sum_{k=0}^{\lfloor \frac{n-1}{2}\rfloor}S(n-k,n-2k)y^k,\]
where $S(n,k)$ are the Stirling numbers of the second kind. It is
worth mentioning that the generating function for the numbers
$g_n(1)$ has been given in OEIS \cite{slo}, that is,
\[\sum_{n\geq 1}g_n(1)x^n=\sum_{n\geq 1}
\frac{x^n}{\prod_{k=1}^n(1-kx^2)}.\] From the generating function of
Stirling numbers of the second kind, it is  straightforward to
deduce that
\[\sum_{n\geq 1}g_n(y) x^n =\sum_{n\geq 1} \frac{x^n}{\prod_{k=1}^n(1-kx^2y)}.
\]

%Using the bijection $\alpha$, we can see that \[g_n(y)=\sum_{k=0}^{{\lfloor \frac{n-1}{2}\rfloor}Q(n-1,k)y^k}\].
%
%%Thus the generating function of partial matchings avoiding neighbor alignments and left nestings is

Below is the recurrence relation of $g_n(y)$
which follows from the recurrence of $S(n,k)$.

%Let $g_n$ denote the coefficient of $x^n$ in
%\[G(x)=\sum_{n\geq 1}
%\frac{x^n}{\prod_{k=1}^n(1-kx^2)}.\]
%The sequence $\{g_n\}_{n\geq 0}$ is list
%as sequence A**** in Sloane [***]. It is known that
%\[g_n =\sum_{k=0}^{\lfloor \frac{n-1}{2}\rfloor}S(n-k,n-2k).\]
%In fact, the above formula can be slightly
%generalized as follows.
\begin{cor}
For $n\geq 3$, $g_n(y)$ has the following recurrence relation
\[g_n(y)=g_{n-1}(y)+(n-2)y\cdot g_{n-2}(y)-2y^2\cdot g'_{n-2}(y).\]
\end{cor}

\section{Neighbor alignments and left, right nestings}

In this section, we obtain
the bivariate generating function for the
number of partial matchings of $[n+k-1]$
with $k$ arcs that avoid
neighbor alignments, left nestings and right nestings. This
generating function turns out to be identical to the
 generating function for the number of self-modified ascent sequences of
 length $n$ with largest element $k-1$
 or $3\bar{1}52\bar{4}$-avoiding permutations of $[n]$ that have $k$
 right-to-left minima
 due to Bousquet-M\'elou, Claesson, Dukes
 and Kitaev \cite{BCDK}.

  Recall that $\R(n,k)$ denotes the set
  of partial matchings of $[n]$ with $k$ arcs that
   avoid neighbor alignments and both left and right nestings.
     We shall give a bijection
      between $\R(2n-k-1,n-k)$ and the set of compositions
of $n-k$ into ${k+1\choose 2}$ components,
 from which we can deduce the
generating function for the numbers $R(n+k-1,k)$. Denote the set of
compositions of $n$ into $k$ components (possibly empty) by
$\C(n,k)$.

\begin{theo}\label{p3-1}
There is a bijection between the set $\R(2n-k-1,n-k)$ and the set
$\C(n-k,{k+1 \choose 2})$.
\end{theo}

\pf Let $M\in \R(2n-k-1,n-k)$ be a partial matching with $2n-k-1$
vertices and $n-k$ arcs containing no left nestings, no right
nestings and no neighbor alignments. We aim to construct a
composition $\beta(M)$ in $\C(n-k,{k+1 \choose 2})$. Clearly, there
are $k-1$ singletons in $M$. These $k-1$ singletons separate the
vertices into $k$ intervals, the first interval is the interval
before the first
 singleton and
the $(i+1)$-st interval is the interval between the $i$-th and
$(i+1)$-st singletons, the $k$-th interval is the interval after the
last singleton.

By the following procedure, we can associate the $i$-th ($1\leq
i\leq k$) interval with an integer composition, that is, a  sequence
$s^{(i)}$ of nonnegative integers of length $k-i+1$.
 For
the origins in the $i$-th interval,
 their corresponding destinations
have $k-i+1$ choices to be in the $i$-th, $(i+1)$-st,$\ldots$, and
the $k$-th interval. If there are $r$ ($r\geq 0$)
 destinations in the $j$-th
($i\leq j\leq k$) interval, then the $(j-i+1)$-st entry
 of the sequence
$s^{(i)}$ is set to be $r$. Since there are $n-k$ destinations in
$M$, thus putting all these $k$ sequences together, we get a
composition
 $s=(s^{(1)},s^{(2)},\ldots,s^{(k)})$
of $n-k$ into $k+(k-1)+\cdots+1={k+1\choose 2}$ components.

Conversely, given a composition of $n-k$ with ${k+1\choose 2}$
components, we may break it into sequences of length $k,k-1,\ldots,1$
respectively, and we denote the $i$-th sequence by
\[ s^{(i)}=(s^{(i)}_i,s^{(i)}_{i+1},
\ldots,s^{(i)}_{k}),\]
where $1\leq i\leq
k.$  Denote the sum of the elements
 of $s^{(i)}$ by $|s^{(i)}|$. We now proceed to
 construct the diagram, or the linear representation
 of a partial matching,   based on the
 sequences  $s^{(1)},s^{(2)},\ldots,s^{(k)}$.  First,
  we draw $k-1$ singletons on a line
to form $k$ intervals such that the first interval is the one
before the first singleton,
 the $(i+1)$-st interval is that between
the $i$-th and $(i+1)$-st singleton, and the $k$-th interval is
 the one after the last
singleton.
Then we need to determine the origins and the
destinations in each interval. We put $|s^{(i)}|$
origins and
$s^{(1)}_i+s^{(2)}_i+\cdots+s^{(i)}_i$ destinations
in the $i$-th
interval, where  all the destinations are located after
all  the origins. So there are $(k-1)+2(n-k)=2n-k-1$
vertices. Next,
we label the vertices from left to right
 by using the numbers
$1,2,\ldots,2n-k-1$.

 Finally, we should match the $n-k$ origins and the $n-k$
 destinations to form $n-k$ arcs. For $i$ from $1$ to $k$, the right-hand endpoints of the
arcs with origins in the $i$-th interval are determined as follows.
As the initial step,  for each $j$ $(i\leq j\leq k)$,  we  choose
the first $s^{(i)}_j$ available destinations (i.e., the destinations
that have not been processed)
 in the $j$-th interval. It is easy to check that
 there are
$s^{(i)}_i+s^{(i)}_{i+1}+\cdots+s^{(i)}_k=|s^{(i)}|$
 destinations that have been chosen so far. Then we match
 these $|s^{(i)}|$ destinations with the $|s^{(i)}|$ origins in the
 $i$-th interval to
form an $|s^{(i)}|$-crossing.
 This construction ensures that there
are neither left nestings nor right nestings.
Furthermore, the
positions of singletons guarantee
 that there are no neighbor
alignments. Therefore we get a partial matching in
$\R(2n-k-1,n-k)$. This implies
that the above mapping $\beta$ is a
bijection, and hence the proof is complete.  \qed

For example, let \[  M=\{\{1,6\},\{2,7\},\{3\},\{4,8\},\{5,14\},
\{9\},\{10,11\},\{12\},\{13,15\}\}\]
 which belongs to $\R(15,6)$. The
three singletons $3,9,12$ separate the vertices into $4$ intervals,
namely, $ \{1,2\},
 \; \{4,5,6,7,8\},\;    \{10,11\},\;
 \{13,14,15\}.$
 For the origins $\{1,2\}$ in the first interval, their
destinations are both in the second interval, so
$s^{(1)}=(0,2,0,0)$. Similarly, we have $s^{(2)}=(1,0,1)$,
$s^{(3)}=(1,0), s^{(4)}=(1)$. So the corresponding
composition of
$6$ into $4+3+2+1=10$ components
 is $(0,2,0,0,1,0,1,1,0,1)$.

\begin{figure}[h,t]
\begin{center}
\begin{picture}(280,220)
\setlength{\unitlength}{0.45mm}

\put(-25,0){$\Longrightarrow $}
\multiput(0,0)(15,0){15}{\circle*{1.5}}
%\put(30,0){\circle{3.5}}\put(120,0){\circle{3.5}}\put(165,0){\circle{3.5}}

\put(-2,-9){\footnotesize1}\put(13,-9){\footnotesize2}
\put(28,-9){\footnotesize3}\put(43,-9){\footnotesize4}\put(58,-9){\footnotesize5}\put(73,-9){\footnotesize6}
\put(88,-9){\footnotesize7}
\put(103,-9){\footnotesize8}\put(118,-9){\footnotesize9}\put(132,-9){\footnotesize10}\put(147,-9){\footnotesize11}
\put(162,-9){\footnotesize12}
\put(177,-9){\footnotesize13}\put(192,-9){\footnotesize14}\put(207,-9){\footnotesize15}

\qbezier[1000](0,0)(37.5,30)(75,0)\qbezier[1000](15,0)(52.5,30)(90,0)
\qbezier[1000](45,0)(75,25)(105,0)\qbezier[1000](60,0)(127.5,45)(195,0)
\qbezier[1000](135,0)(142.5,20)(150,0)\qbezier[1000](180,0)(195,20)(210,0)

%%%%%%%%%%%%%%%%%%%%%%%%%%%%%%%%%%%%%%%%%%%%%%%%%%%%%%%%%%%%%%%%%%%%%
%\put(0,0){$\line(1,1){8}$} \put(0,0){$\line(-1,1){8}$}
\put(-25,50){$\Longrightarrow $}
\multiput(0,50)(15,0){15}{\circle*{1.5}}
%\put(30,50){\circle{3.5}}\put(120,50){\circle{3.5}}\put(165,50){\circle{3.5}}

\put(-2,41){\footnotesize1}\put(13,41){\footnotesize2}
\put(28,41){\footnotesize3}\put(43,41){\footnotesize4}\put(58,41){\footnotesize5}\put(73,41){\footnotesize6}\put(88,41){\footnotesize7}
\put(103,41){\footnotesize8}\put(118,41){\footnotesize9}\put(132,41){\footnotesize10}\put(147,41){\footnotesize11}\put(162,41){\footnotesize12}
\put(177,41){\footnotesize13}\put(192,41){\footnotesize14}\put(207,41){\footnotesize15}

\put(135,50){$\line(1,1){6}$} \put(150,50){$\line(-1,1){6}$}
\put(180,50){$\line(1,1){6}$} \put(210,50){$\line(-1,1){6}$}
\qbezier[1000](0,50)(37.5,85)(75,50)\qbezier[1000](15,50)(52.5,85)(90,50)
\qbezier[1000](45,50)(75,80)(105,50)\qbezier[1000](60,50)(122.5,95)(195,50)

%%%%%%%%%%%%%%%%%%%%%%%%%%%%%%%%%%%%%%%%%%%%%%%%%%%%%%%%%%%%%%%%%%%%%
\put(-25,100){$\Longrightarrow $}
\multiput(0,100)(15,0){15}{\circle*{1.5}}
%\put(30,100){\circle{3.5}}\put(120,100){\circle{3.5}}\put(165,100){\circle{3.5}}

\put(-2,91){\footnotesize1}\put(13,91){\footnotesize2}
\put(28,91){\footnotesize3}\put(43,91){\footnotesize4}\put(58,91){\footnotesize5}\put(73,91){\footnotesize6}\put(88,91){\footnotesize7}
\put(103,91){\footnotesize8}\put(118,91){\footnotesize9}\put(132,91){\footnotesize10}\put(147,91){\footnotesize11}\put(162,91){\footnotesize12}
\put(177,91){\footnotesize13}\put(192,91){\footnotesize14}\put(207,91){\footnotesize15}

\put(45,100){$\line(1,1){6}$}\put(60,100){$\line(1,1){6}$}
\put(105,100){$\line(-1,1){6}$}\put(135,100){$\line(1,1){6}$}
\put(150,100){$\line(-1,1){6}$}\put(180,100){$\line(1,1){6}$}
\put(195,100){$\line(-1,1){6}$}\put(210,100){$\line(-1,1){6}$}

\qbezier[1000](0,100)(37.5,135)(75,100)
\qbezier[1000](15,100)(52.5,135)(90,100)

%%%%%%%%%%%%%%%%%%%%%%%%%%%%%%%%%%%%%%%%%%%%%%%%%%%%%%%%%%%%%%%%%%%%%
\multiput(0,150)(15,0){15}{\circle*{1.5}}
%\put(30,150){\circle{3.5}}\put(120,150){\circle{3.5}}
%\put(165,150){\circle{3.5}}

\put(-2,141){\footnotesize1}\put(13,141){\footnotesize2}
\put(28,141){\footnotesize3}\put(43,141){\footnotesize4}\put(58,141){\footnotesize5}
\put(73,141){\footnotesize6}\put(88,141){\footnotesize7}
\put(103,141){\footnotesize8}\put(118,141){\footnotesize9}\put(132,141){\footnotesize10}
\put(147,141){\footnotesize11}\put(162,141){\footnotesize12}
\put(177,141){\footnotesize13}\put(192,141){\footnotesize14}\put(207,141){\footnotesize15}

\put(0,150){$\line(1,1){6}$}\put(15,150){$\line(1,1){6}$}
\put(45,150){$\line(1,1){6}$}\put(60,150){$\line(1,1){6}$}
\put(75,150){$\line(-1,1){6}$}\put(90,150){$\line(-1,1){6}$}
\put(105,150){$\line(-1,1){6}$}\put(135,150){$\line(1,1){6}$}
\put(150,150){$\line(-1,1){6}$} \put(180,150){$\line(1,1){6}$}
\put(195,150){$\line(-1,1){6}$}\put(210,150){$\line(-1,1){6}$}
\end{picture}
\vspace*{5mm} \caption{The bijection $\beta$.}\label{fig5}
\end{center}
\end{figure}
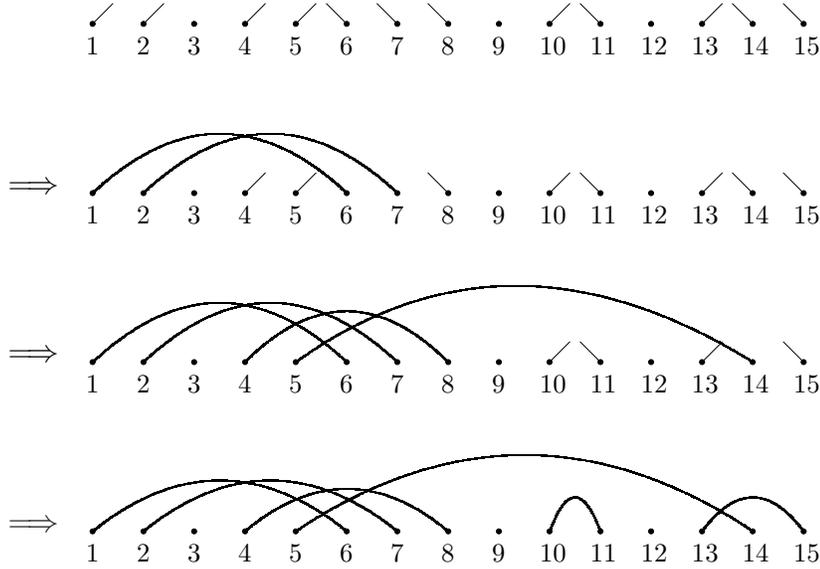

Conversely, given the composition $(0,2,0,0,1,0,1,1,0,1)$,
we split
it into $4$ sequences,
$s^{(1)}=(s^{(1)}_1,s^{(1)}_2,s^{(1)}_3,s^{(1)}_4)=(0,2,0,0),
 s^{(2)}=(s^{(2)}_2,s^{(2)}_3,s^{(2)}_4)
=(1,0,1), s^{(3)}=(s^{(3)}_3,s^{(3)}_4)=(1,0),
s^{(4)}=(s^{(4)}_4)=(1)$. The  construction
of the corresponding partial matching is illustrated in Figure \ref{fig5}.

{\noindent \it Proof of Theorem \ref{t3}.}   Note that the
coefficient of $x^n$ in \[\sum_{n\geq 1} \frac{x^n}{(1-xy)^{{n+1
\choose 2}}}\]
 equals \[\sum_{k=0}^n{{k \choose 2}+n-1 \choose
n-k}y^{n-k},\] which equals the number of compositions of $n-k$ into
${k+1 \choose 2}$ components. Thus the result follows from Theorem
\ref{p3-1}. \qed

\section{Partitions with no right nestings}

The objective of this section is to construct a bijection between
the set $\T(n-k+1,k)$ of partitions of $[n-k+1]$ with $k$ arcs but
with no right nestings and the set $\R(n,k)$ of partial matchings of
$[n]$ with $k$ arcs but with no left nestings, right nestings and
neighbor alignments. In fact, we only need to establish a
correspondence between a sequence of $n-2k+1$ compositions and the
set $\T(n-k+1,k)$. Combining the bijection $\beta$ in Section 4 from
compositions to partial matchings without left,
 right nestings and
neighbor alignments, we
 obtain the desired bijection between the set $\R(n,k)$ and the set
$\T(n-k+1,k)$. Since in the previous section we have computed the
generating function for the numbers $R(n,k)$,  we are led to the
generating function for $T(n,k)$ as stated in Theorem \ref{t4}.

\begin{theo}\label{p3-2}
There exists a bijection between the set $\R(n,k)$ and the set
$\T(n-k+1,k)$. Moreover, this bijection transforms the number of
left crossings  of a partial matching into the number of transients of
a partition.
\end{theo}

\pf  Let $M\in\R(n,k)$, that is, $M$ is a partial matching of
$[n]$ with $k$ arcs but with no left nestings, right nestings and
neighbor alignments. We wish to construct a partition
$\gamma(M)\in\T(n-k+1,k)$. The construction consists of two steps.
The first step is to associate the partial matching $M\in\R(n,k)$
with a sequence  of $n-2k+1$ compositions. The second step is to
obtain the desired partition $\gamma(M)\in\T(n-k+1,k)$ from those
compositions.

Given a partial matching $M\in\R(n,k)$, intuitively the $n-2k$
singletons of $M$ break the vertices of $M$
 into $n-2k+1$ intervals since the vertices are assumed to
 be arranged in increasing order on the horizontal line.
 As in the
construction of the bijection $\beta$
between compositions and partial matchings, we
associate the $i$-th interval with a composition
$s^{(i)}=(s^{(i)}_i,\ldots,s^{(i)}_{n-2k+1})$,
where $s^{(i)}_j
(i\leq j\leq n-2k+1)$ is the number of arcs with origins in the
$i$-interval and destinations in the $j$-th interval.

The procedure to generate the partition $\gamma(M)$  can be
described  as follows. We start with $n-2k+1$
 empty intervals by putting
down $n-2k$ singletons on a line. Then we determine the left-hand
and right-hand endpoints of every interval so that all the arcs
are consequently determined by the endpoints.

To reach this goal, we define a $k$-path  to be a sequence of $k$
arcs of the form $(i_1,i_2),(i_2,i_3),$ $\ldots, (i_k,i_{k+1})$,
where $i_1<i_2<\cdots <i_{k+1}$. For $i$ from $1$ to $n-2k+1$, we
construct an $|s^{(i)}|$-path $(i_1,i_2),(i_2,i_3),\ldots,$
$(i_{|s^{(i)}|},i_{|s^{(i)}|+1})$ via the following steps.

\noindent Step 1. We put the origin $i_1$ of this path immediately
 before the leftmost right-hand endpoint
  that has been constructed in the $i$-th interval. If there are no
  right-hand endpoints in the $i$-th interval, just put $i_1$ before
  the $i$-th singleton.

\noindent Step 2.  According to the composition $s^{(i)}$,
 we determine the
positions of the right-hand endpoints $i_2,\ldots,i_{|s^{(i)}|+1}$
of this path. More precisely, we assign $s^{(i)}_j (i\leq j\leq
n-2k+1)$ right-hand endpoints to the $j$-th $(i\leq j\leq n-2k+1)$
interval.

The $|s^{(i)}|$-path is constructed by inserting the arcs
$(i_1,i_2),(i_2,i_3),\ldots,(i_{|s^{(i)}|},i_{|s^{(i)}|+1})$ one by
one. Precisely, by inserting an arc to an interval we mean inserting
the right-hand endpoint of this arc to this interval. We claim that
the positions of the right-hand endpoints
$i_2,\ldots,i_{|s^{(i)}|+1}$ of the $|s^{(i)}|$-path in each
interval are uniquely determined subject to the
constraint that no  right nestings are allowed.
 To prove this
claim, we show that for each arc $(i_s,i_{s+1})$ ($1\leq s\leq
|s^{(i)}|$),
 there is one and only
one position to insert the right-hand endpoint $i_{s+1}$.

Suppose that we wish to insert the arc $e=(i_s,i_{s+1})$ to the
$j$-th interval, where the left-hand endpoint $i_s$ of $e$ is
determined already. We proceed to determine the position of
$i_{s+1}$. If there are no right-hand endpoints to the right of
$i_s$ in this interval, then we insert $i_{s+1}$ immediately before
the $j$-th singleton. Otherwise, we assume that there are $t$
right-hand endpoints $r_1,r_2,\ldots,r_{t-1},r_t$ to the right of
 $i_s$ in the
$j$-th interval. As will be seen, there is a unique position
 to insert  $e$ to the $j$-th
interval such that no right nestings will be formed.

The strategy of inserting $e$ can be easily described as follows. We
begin with the position immediately to the left of $r_1$.
 If $i_{s+1}$ can be placed in this position without
  causing any right
nestings, then this is the position
we are looking for. Otherwise, we   consider the
position immediately before $r_2$ as the second candidate.

Like the case for $r_1$, if putting $i_{s+1}$  immediately before
$r_2$
 does not cause any
right nestings, then it is the desired choice.
 Otherwise, we
consider the position immediately before $r_3$
 as the third
candidate. Then we continue this process until we
 find a position such that
after inserting $e$ creates no right nestings.

To see that the above process will terminate
at some point, we assume that  $i_{s+1}$
cannot be inserted immediately before $r_i$, and we
assume that
 inserting $i_{s+1}$ immediately after $r_i$ also
 yields  a right nesting. Then this right nesting
 caused by the insertion of $i_{s+1}$ immediately after $r_i$
  must be formed
  by the arc $e$ and the arc whose
   right-hand endpoint is  immediately after $r_i$.
   This means that there is a right-hand endpoint after $r_i$.
Since the
 number of right-hand endpoints in every interval is finite,
  we conclude  that
 there always exits a  position such that inserting
  $i_{s+1}$ does not cause any right nestings.

It is still necessary to  show that there
is a unique choice for the position
of $i_{s+1}$. Assume
 that we have found a position immediately before
the vertex $r_{i_1}$ such that
the insertion of $e$ does not cause
any right nestings. It can be shown that all the
positions to the right of $r_{i_1}$ cannot be chosen
for the insertion of $e$. Otherwise, suppose that the
position immediately after the vertex $r_{i_2}$ is a
feasible choice, where $r_{i_1} < r_{i_2}$.

We now proceed to find a right nesting that implies
a contradiction. If
$i_{s+1}$ can be inserted immediately before $r_{i_1}$,
 then the arc $e$ and
the arc $e_1=(l_1,r_{i_1})$  form a crossing, that is,
$i_s<l_1$; On the other hand, if $i_{s+1}$ can be inserted
immediately after $r_{i_2}$, then  $e$ and
$e_2=(l_2,r_{i_2})$   form a crossing as well, that is,
$l_2<i_s$. This implies that $l_1>l_2$.
So  the   arcs $e_1$ and $e_2$
form a nesting.

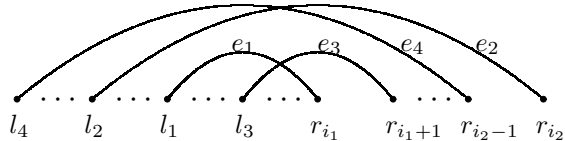
\begin{figure}[h,t]
\begin{center}
\begin{picture}(180,50)
\setlength{\unitlength}{0.5mm}

\multiput(0,0)(20,0){8}{\circle*{1.5}}
\multiput(6,-2)(20,0){4}{$\cdots$} \put(106,-2){$\cdots$}
\put(-2,-9){\footnotesize$l_4$}\put(18,-9){\footnotesize$l_2$}
\put(38,-9){\footnotesize$l_1$}\put(58,-9){\footnotesize$l_3$}
\put(78,-9){\footnotesize$r_{i_1}$}\put(98,-9){\footnotesize$r_{i_1+1}$}
\put(118,-9){\footnotesize$r_{i_2-1}$}\put(138,-9){\footnotesize$r_{i_2}$}

\qbezier[1000](40,0)(60,25)(80,0) \qbezier[1000](20,0)(80,50)(140,0)
\qbezier[1000](60,0)(80,25)(100,0)\qbezier[1000](0,0)(60,50)(120,0)
\put(57,13){\footnotesize$e_1$}\put(122,13){\footnotesize$e_2$}
\put(80,13){\footnotesize$e_3$}\put(102,13){\footnotesize$e_4$}

\end{picture}
\vspace*{.5cm}
 \caption{The uniqueness of inserting an arc.}\label{fig6}
\end{center}
\end{figure}

In order to find a right nesting,
we consider the distance between
$r_{i_1}$ and $r_{i_2}$.
 If $r_{i_1}+1=r_{i_2}$, then the two arcs $e_1$ and $e_2$ form a right
nesting.
 If $r_{i_1}+2=r_{i_2}$, namely,
 there is a vertex $r_{i_1+1}$ between
$r_{i_1}$ and $r_{i_2}$, then the arc with right-hand endpoint
$r_{i_1+1}$ forms a right nesting with the arc $e_1$
 or $e_2$.
We now consider the case that
 there are more than one vertices between $r_{i_1}$
 and $r_{i_2}$.
Since in every step of the inserting process no right
nestings are formed,
the left-hand endpoint $l_3$ of the arc
$e_3=(l_3,r_{i_1+1})$ must be  to the right of $l_1$, and
the
position of the left-hand endpoint $l_4$ of the arc
$e_4=(l_4,r_{i_2-1})$ must be to the left of $l_2$.
Thus we deduce that  $e_3$
and $e_4$ form a nesting as well, and the distance
between the right-hand endpoints of $e_3$ and
 $e_4$ has decreased by two compared with the
 distance between the right-hand endpoints of $e_1$ and
 $e_2$. See Figure \ref{fig6} for an illustration.

Iterating the above process by
 checking the distance between the
point $r_{i_1+1}$ and the point $r_{i_2-1}$,
 we can always find
a right nesting. This is a contradiction,
which means that  there is a unique choice for
the insertion of $e$ without causing right nestings.

By the above uniqueness property,
we can insert the arcs
$(i_1,i_2),(i_2,i_3),\ldots,(i_{|s^{(i)}|},i_{|s^{(i)}|+1})$
one by
one to construct a unique $|s^{(i)}|$-path.
 After  $n-2k+1$
steps, we get a partition with no right nestings. Finally, delete
every singleton that is immediately to the left of an origin, except
for the first origin. Examining  the number of points as in the
construction of the bijection $\alpha$ in Section 3,
 we are  led
to a partition $\gamma(M)$ on $[n-k+1]$ without right nestings.

Conversely, given a  partition
$\pi\in \T(n+1-k,k)$ with $k$ arcs
that has no right nestings,
 we wish to construct a partial matching in
$\R(n,k)$. Clearly, we should add $k-1$ vertices into $\pi$.
 First of all, we
add a vertex before each origin except the first
origin. At this point, the number of added vertices
is the number of non-singleton blocks of $\pi$ minus one.
Assume that the new
partition $\pi'$ has $m$ singletons
 which split the vertex set
into $m+1$ intervals.

In each interval, there is at most one origin.
Assume that the origin in the $i$-th interval
is the origin of an
$r$-path, then we associate the $i$-th interval with
 a composition
$t^{(i)}=(t^{(i)}_i,\ldots,t^{(i)}_{m+1})$ of the integer
$r$,  where $t^{(i)}_j (i\leq
j\leq m+1)$ is the number of right-hand endpoints
of this $r$-path
in the $j$-th interval.
From these $m+1$ compositions, using the map
$\beta$ in Section 4 from compositions to
partial matchings without
left, right nestings and neighbor alignments,
 we  obtain a
partial matching of $\R(n,k)$ with $k$ arcs.
 It is easily seen  that the number of left crossings of
$M$ equals the number of transients of $\gamma(M)$.
 This completes
the proof. \qed

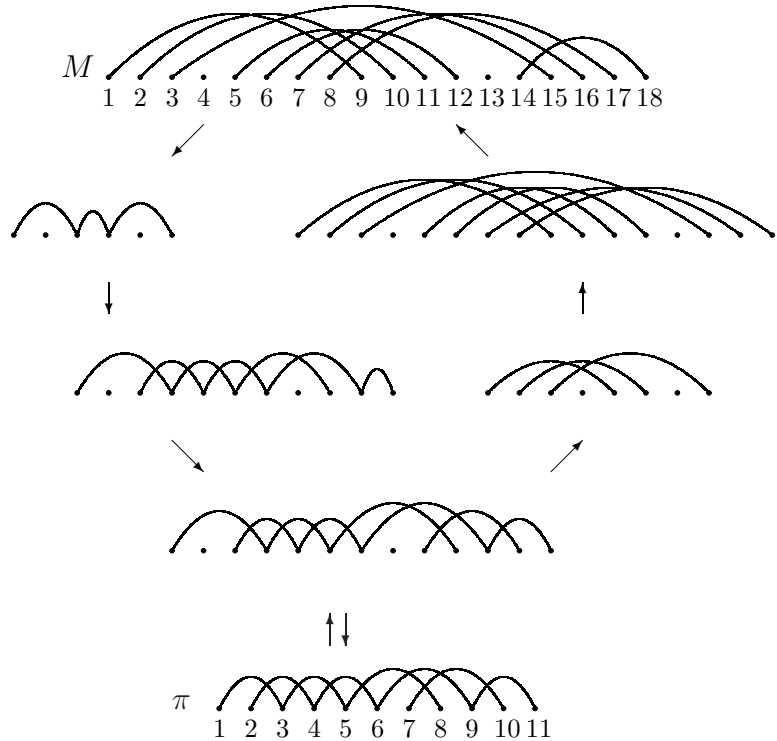
\begin{figure}[h,t]
\begin{center}
\begin{picture}(180,280)
\setlength{\unitlength}{0.42mm}
%%%%%%%%%%%%%%%%%%%%%%%%%%%%%%%%%%%%%%%%%%%%%%%%%%%%%%%%%%%%%%%%%%%%%

\put(50,20){\vector(0,1){10}}\put(55,30){\vector(0,-1){10}}

\put(0,0){$\pi$} \multiput(15,0)(10,0){11}{\circle*{1.5}}

\put(13,-9){\footnotesize1}\put(23,-9){\footnotesize2}\put(33,-9){\footnotesize3}\put(43,-9){\footnotesize4}\put(53,-9){\footnotesize5}\put(63,-9){\footnotesize6}
\put(73,-9){\footnotesize7}\put(83,-9){\footnotesize8}\put(93,-9){\footnotesize9}\put(102,-9){\footnotesize10}\put(112,-9){\footnotesize11}

\qbezier[1000](15,0)(25,20)(35,0) \qbezier[1000](25,0)(35,20)(45,0)
\qbezier[1000](45,0)(55,20)(65,0)\qbezier[1000](95,0)(80,25)(65,0)

 \qbezier[1000](95,0)(105,20)(115,0)
\qbezier[1000](55,0)(45,20)(35,0)\qbezier[1000](55,0)(70,25)(85,0)\qbezier[1000](75,0)(90,25)(105,0)
%%%%%%%%%%%%%%%%%%%%%%%%%%%%%%%%%%%%%%%%%%%%%%%%%%%%%%%%%%%%%%%%%%%%%
 \multiput(0,50)(10,0){13}{\circle*{1.5}}

%\put(-2,41){1}\put(8,41){2}\put(18,41){3}\put(28,41){4}\put(38,41){5}
%\put(48,41){6}\put(58,41){7}\put(68,41){8}\put(78,41){9}\put(87,41){10}
%\put(97,41){11}\put(107,41){12}\put(117,41){13}

\qbezier[1000](0,50)(15,75)(30,50)\qbezier[1000](50,50)(40,70)(30,50)
\qbezier[1000](50,50)(70,80)(90,50)\qbezier[1000](20,50)(30,70)(40,50)
\qbezier[1000](60,50)(50,70)(40,50)\qbezier[1000](60,50)(80,80)(100,50)

\qbezier[1000](120,50)(110,70)(100,50)
%\color{red}
\qbezier[1000](110,50)(95,75)(80,50)
%%%%%%%%%%%%%%%%%%%%%%%%%%%%%%%%%%%%%%%%%%%%%%%%%%%%%%%%%%%%%%%%%%%%%
%%%%%%%%%%%
%\color{black}
\put(0,85){\vector(1,-1){10}}\put(120,75){\vector(1,1){10}}

\multiput(-30,100)(10,0){11}{\circle*{1.5}}

\qbezier[1000](-30,100)(-15,125)(0,100)

 \qbezier[1000](0,100)(10,120)(20,100)
\qbezier[1000](20,100)(35,125)(50,100)
%\color{red}
\qbezier[1000](-10,100)(0,120)(10,100)\qbezier[1000](10,100)(20,120)(30,100)
\qbezier[1000](60,100)(65,115)(70,100)\qbezier[1000](60,100)(45,125)(30,100)
%\color{black}

 \multiput(100,100)(10,0){8}{\circle*{1.5}}
%\put(98,91){1}\put(108,91){2}\put(118,91){3}\put(128,91){4}\put(138,91){5}
%\put(148,91){6}\put(158,91){7}\put(168,91){8}
\qbezier[1000](100,100)(120,120)(140,100)\qbezier[1000](110,100)(130,120)(150,100)
\qbezier[1000](120,100)(145,125)(170,100)

%%%%%%%%%%%%%%%%%%%%%%%%%%%%%%%%%%%%%%%%%%%%%%%%%%%%%%%%%%%%%%%%%%%%%
\put(10,185){\vector(-1,-1){10}}\put(100,175){\vector(-1,1){10}}
\put(-20,135){\vector(0,-1){10}}\put(130,125){\vector(0,1){10}}
\multiput(-50,150)(10,0){6}{\circle*{1.5}}
%\put(-52,141){1}\put(-42,141){2}\put(-32,141){3}\put(-22,141){4}\put(-12,141){5}
%\put(-2,141){6}
\qbezier[1000](-50,150)(-40,170)(-30,150)\qbezier[1000](-20,150)(-25,165)(-30,150)
\qbezier[1000](-20,150)(-10,170)(0,150)

\multiput(40,150)(10,0){16}{\circle*{1.5}}
%\put(38,141){1}\put(48,141){2}\put(58,141){3}\put(68,141){4}\put(78,141){5}
%\put(88,141){6}\put(98,141){7}\put(108,141){8}\put(118,141){9}\put(127,141){10}
%\put(137,141){11}\put(147,141){12}\put(157,141){13}
%\put(167,141){14}\put(177,141){15}\put(187,141){16}

\qbezier[1000](40,150)(80,185)(120,150)\qbezier[1000](50,150)(90,185)(130,150)
\qbezier[1000](60,150)(115,190)(170,150)

%\color{red}
\qbezier[1000](80,150)(110,180)(140,150)
\qbezier[1000](90,150)(120,180)(150,150)
\qbezier[1000](100,150)(140,180)(180,150)\qbezier[1000](110,150)(150,180)(190,150)
%%%%%%%%%%%%%%%%%%%%%%%%%%%%%%%%%%%%%%%%%%%%%%%%%%%%%%%%%%%%%%%%%%%%%
\put(-35,200){$M$} \multiput(-20,200)(10,0){18}{\circle*{1.5}}

\put(-22,191){\footnotesize1}\put(-12,191){\footnotesize2}\put(-2,191){\footnotesize3}\put(8,191){\footnotesize4}\put(18,191){\footnotesize5}
\put(28,191){\footnotesize6}\put(38,191){\footnotesize7}\put(48,191){\footnotesize8}\put(58,191){\footnotesize9}\put(67,191){\footnotesize10}
\put(77,191){\footnotesize11}\put(87,191){\footnotesize12}\put(97,191){\footnotesize13}
\put(107,191){\footnotesize14}\put(117,191){\footnotesize15}\put(127,191){\footnotesize16}
\put(137,191){\footnotesize17}\put(147,191){\footnotesize18}

\qbezier[1000](-20,200)(20,240)(60,200)
\qbezier[1000](-10,200)(30,240)(70,200)
\qbezier[1000](0,200)(60,245)(120,200)
\qbezier[1000](20,200)(50,230)(80,200)
\qbezier[1000](30,200)(60,230)(90,200)
\qbezier[1000](40,200)(80,240)(130,200)
\qbezier[1000](50,200)(90,240)(140,200)
\qbezier[1000](110,200)(130,225)(150,200)

%\put(-40,230){$s^1=(0,2,1)$} \put(30,230){$s^2=(2,2)$}
%\put(120,230){$s^3=(1)$}
%
%\put(10,190){\line(0,1){60}}\put(100,190){\line(0,1){60}}
\end{picture}
\vspace*{5mm} \caption{The bijection $\gamma$.
}\label{fig7}
\end{center}
\end{figure}

Figure \ref{fig7} gives an example of a
 partial matching $M$  without left,
 right nestings and neighbor
alignments. It also illustrates the process to construct a
partition $\gamma(M)$ without right nestings.
There are two singletons in $M$ which create
three intervals. The first interval is associated with the
composition $s^{(1)}=(0,2,1)$, which is transformed
into a $3$-path
of $\gamma(M)$.
The second interval is associated with the
composition $s^{(2)}=(2,2)$, which is transformed into
 a $4$-path of
$\gamma(M)$. The third interval is associated
with the composition
$s^{(3)}=(1)$, which is transformed into a
$1$-path of $\gamma(M)$.

To conclude, we remark that in general the number
of  partitions of $[n]$ avoiding right
crossings is not equal to the number of  partitions
of $[n]$ avoiding right nestings.
It would be interesting to find
the generating function for the
number of  partitions of $[n]$
without right crossings.

\vspace{.2cm} \noindent{\bf Acknowledgments.} This work was
supported by  the 973 Project, the PCSIRT Project of the Ministry of
Education, and the National
Science Foundation of China.

%-------------------------------------------------------------


\begin{thebibliography}{99}

\bibitem{BCDK}
M. Bousquet-M\'elou, A. Claesson, M. Dukes and S. Kitaev,
$(2+2)$-free posets, ascent sequences and pattern avoiding
permutations, J. Combin. Theory Ser. A, 117 (2010), 884--909.


\bibitem{Chen05}
W.Y.C. Chen, E.Y.P. Deng and R.R.X. Du, Reduction of $m$-regular
noncrossing partitions, European J. Combin., 26 (2005), 237--243.

\bibitem{Chen06}
W.Y.C. Chen, T. Mansour and S.H.F. Yan, Matchings avoiding
partial patterns, Electron. J. Combin., 13 (2006), R112.


\bibitem{Chen07}
W.Y.C. Chen, E.Y.P. Deng, R.R.X. Du, R.P. Stanley and C.H.
Yan, Crossings and nestings of matchings and partitions, Trans.
Amer. Math. Soc., 359 (2007), 1555--1575.


\bibitem{CL}
A. Claesson and S. Linusson, $n!$ matchings, $n!$ posets, Proc.
Amer. Math. Soc., to appear.



\bibitem{CDK}
A. Claesson, M. Dukes and S. Kitaev, A direct encoding of
Stoimenow's matchings as ascent sequences,
arXiv:math.CO/0910.1619.

\bibitem{M}
M. de Sainte-Catherine, Couplages et Pfaffiens en combinatoire,
physique et informatique, Ph.D. Thesis, University of Bordeaux I,
1983.

\bibitem{JLMY}
V. Jelinek, N.Y. Li, T. Mansour and S.H.F. Yan, Matchings
avoiding partial patterns and lattice paths, Electron. J. Combin., 13
(2006), R89.


\bibitem{Kla1}
M. Klazar, On $abab$-free and $abba$-free set partitions, European J.
Combin., 17 (1996), 53--68.

\bibitem{Kla98}
M. Klazar, On trees and noncrossing partitions, Discrete Appl. Math., 82 (1998), 263--269.



\bibitem{Kla2}
M. Klazar,  Bell numbers, their relatives, and algebraic differential
equations, J. Combin. Theory Ser. A, 102 (2003), 63--87.

\bibitem{Kla3}
M. Klazar, Non-P-recursiveness of   numbers of matchings or linear chord
diagrams with many crossings, Adv. Appl. Math., 30 (2003), 126--136.


\bibitem{slo}
N.J.A. Sloane, The On-Line Encyclopedia of Integer Sequences,
http://www.research.att.com/snjas/sequences.


\bibitem{Sta}
R.P. Stanley, Enumerative Combinatorics, Vol. 1, Cambridge
University Press, Cambridge, UK, 1996.

\bibitem{Sto}
A. Stoimenow, Enumeration of chord diagrams and an upper bound for
Vassiliev invariants, J. Knot Theory Ramifications, 7 (1998),
93--114.

\bibitem{Stein}
P.R. Stein, On a class of linked diagrams, I. Enumerations, J.
Combin. Theory Ser. A, 24 (1978), 357--366.


\end{thebibliography}
\end{document}